\documentclass{amsart}

\usepackage{amscd,amssymb}
\usepackage[all]{xy}

\def\Z{{\mathbb Z}}
\def\Q{{\mathbb Q}}
\def\R{{\mathbb R}}
\def\C{{\mathbb C}}
\def\F{{\mathbb F}}

\def\A{{\mathcal A}}
\def\B{{\mathcal B}}
\def\H{{\mathcal H}}
\def\J{{\mathcal J}}

\def\M{{\mathcal M}}
\def\S{{\mathcal S}}
\def\T{{\mathcal T}}
\def\X{{\mathcal X}}
\def\Y{{\mathcal Y}}

\def\cQ{{\mathcal Q}}

\def\h{{\mathfrak h}}

\def\ubar{\overline{u}}

\def\G{{\Gamma}}
\def\bgamma{{\vec\gamma}}

\def\Mbar{\overline{\M}}

\def\dot{{\bullet}}
\def\bs{\backslash}

\def\red{\mathrm{red}}
\def\sing{\mathrm{sing}}
\def\sep{\mathrm{sep}}
\def\CW{\mathrm{CW}}

\def\liminj#1{\lim_{\stackrel{\longrightarrow}{#1}}}

\renewcommand\Im{\operatorname{Im}}

\newcommand\im{\operatorname{im}}               

\newcommand\Hom{\operatorname{Hom}}
\newcommand\Aut{\operatorname{Aut}}
\newcommand\Diff{\operatorname{Diff}}
\newcommand\Jac{\operatorname{Jac}}
\newcommand\SL{\operatorname{SL}}
\newcommand\Sp{\operatorname{Sp}}

\newtheorem{theorem}{Theorem}
\newtheorem{lemma}[theorem]{Lemma}
\newtheorem{proposition}[theorem]{Proposition}
\newtheorem{corollary}[theorem]{Corollary}

\theoremstyle{definition}

\newtheorem{conjecture}{Conjecture}
\newtheorem{problem}{Problem} 

\theoremstyle{remark}
\newtheorem{remark}[theorem]{Remark}

\begin{document}

\title{Finiteness and Torelli Spaces}

\author{Richard Hain}
\address{Department of Mathematics\\ Duke University\\
Durham, NC 27708-0320}
\email{hain@math.duke.edu}


\date{\today}

\thanks{Supported in part by grants from the National Science Foundation.}


\maketitle

Torelli space $\T_g$ ($g\ge 2$) is the quotient of Teichm\"uller space by the
Torelli group $T_g$. It is the moduli space of compact, smooth genus $g$ curves
$C$ together with a symplectic basis of $H_1(C;\Z)$ and is a model of the
classifying space of $T_g$. Mess, in his thesis \cite{mess}, proved that $\T_2$
has the homotopy type of a bouquet of a countable number of circles. Johnson
and Millson (cf.\ \cite{mess}) pointed out that a similar argument shows that
$H_3(\T_3)$ is of infinite rank. Akita \cite{akita} used an indirect argument
to prove that $\T_g$ does not have the homotopy type of a finite complex for
(almost) all $g\ge 2$. However, the infinite topology of $\T_g$ is not well
understood. The results of Mess and Johnson-Millson are the only ones I know of
that explicitly describe some infinite topology of any Torelli space. Moreover,
although Johnson \cite{johnson:fg} proved that $T_g$ is finitely generated when
$g\ge 3$, there is not one $g\ge 3$ for which it is known whether $T_g$ is
finitely presented or not.

The goal of this note is to present a suite of problems designed to probe the
infinite topology of Torelli spaces in all genera. These are presented in the
fourth section of the paper. The second and third sections present background
material needed in the discussion of the problems.

To create a context for these problems, we first review the arguments of Mess
and Johnson-Millson. As explained in Section~\ref{geography}, Torelli space in
genus 2 is the complement of a countable number of disjoint smooth divisors
$D_\alpha$ (i.e., codimension 1 complex subvarieties) in $\h_2$, the Siegel
upper half plane of rank 2. More precisely,
$$
\T_2 = \h_2 - \bigcup_{\phi \in \Sp_2(\Z)} \phi (\h_1\times \h_1)
$$
where $\Sp_g(R)$ denotes the automorphisms of $R^{2g}$ that fix the standard
symplectic inner product.\footnote{Note that $\Sp_1 = \SL_2$.} One can now
argue, as Mess did \cite{mess}, that $\T_2$ is homotopy equivalent to a
countable wedge of circles. Such an explicit description of $\T_g$ is unlikely
in higher genus due to the increasing complexity of the image of the period
mapping and the lack of an explicit description of its closure. So it should be
fruitful to ponder the source of the infinite rank homology without appeal to
this explicit description of $\T_2$. To this end, consider the exact sequence
$$
\cdots
\to H_c^4(\h_2) \to H_c^4(D) \to H_c^5(\T_2) \to H_c^5(\h_2) \to
\cdots
$$
where $D=\cup D_\alpha$ is a locally finite union of smooth divisors indexed by
$\Sp_2(\Z)$ mod the stabilizer $S_2\ltimes \Sp_1(\Z)^2$ of $\h_1\times
\h_1$.\footnote{The symmetric group on 2 letters, $S_2$, acts on $\h_1\times
\h_1$ by swapping the two factors.} Since $\h_2$ is a contractible complex
manifold of real dimension 6,
$$
H_c^k(\h_2) = H_{6-k}(\h_2) = 0 \text{ when } k < 6.
$$
Thus
$$
H_1(\T_2) = H^5_c(\T_2) \cong H_c^4(D) \cong
\bigoplus_{\alpha \in \Sp_2(\Z)/S_2\ltimes \Sp_1(\Z)^2} \Z[D_\alpha],
$$
a free abelian group of countable rank.

The situation is genus 3 is similar, but slightly more complicated. In this
case, the period mapping $\T_3 \to \h_3$ is two-to-one, branched along the
locus of hyperelliptic curves. The image is the complement $\J_3$ of a
countable union $Z$ of submanifolds $Z_\alpha$ of complex codimension 2:
$$
\J_3 = \h_3 - Z =
\h_3 - \bigcup_{\phi \in \Sp_3(\Z)/\Sp_1(\Z)\times\Sp_2(\Z)}
\phi(\h_1 \times \h_2);
$$
it is the set of framed jacobians of smooth genus 3 curves. By elementary
topology,
$$
H_\dot(\J_3;\Q) \cong H_\dot(\T_3;\Q)^{\Z/2\Z}
$$
where $\Z/2\Z$ is the group of automorphisms of the ramified covering $\T_3 \to
\J_3$. If $H_k(\J_3;\Q)$ is infinite dimensional, then so is $H_k(\T_3;\Q)$. As
in the genus 2, we have an exact sequence
$$
0 = H^8_c(\h_3) \to H^8_c(Z) \to H^9_c(\J_3) \to H^9_c(\h_3) = 0
$$
so that
$$
H_3(\J_3) = H^9_c(\J_3) \cong H^8_c(Z) =
\bigoplus_{\alpha \in \Sp_3(\Z)/\Sp_1(\Z)\times\Sp_2(\Z)} \Z[Z_\alpha].
$$

Three important ingredients in these arguments are:
\begin{itemize}

\item $\T_g$ or a space very closely related to it (e.g., $\J_g$) is the
complement $X-Z$ in a manifold $X$ of a countable union $Z = \cup Z_\alpha$ of
smooth subvarieties;

\item The topology of the manifold $X$ is very simple --- in the cases above,
$X=\h_g$, which is contractible;

\item the ``topology at infinity'' of $X$ is very simple --- in the cases above
the boundary of $X = \h_g$ is a sphere;

\item the topology of $X$, $Z$ and $\T_g$ are related by a Gysin sequence of
compactly supported cohomology, and the relationship between the compactly
supported cohomology and ordinary cohomology is entwined with the ``topology at
infinity'' of $X$.

\end{itemize}

What is special in genus 2 and 3 is that the closure $\J_g^c$ of the image of
the period mapping $\T_g \to \h_g$ is very simple --- it is all of $\h_g$, a
topological ball, which is contractible and has very simple topology at
infinity. This fails to be true in genus 4, where $\J_g^c$ is a singular
subvariety of complex codimension 1. Here, even though there is an explicit
equation for the $\J_4^c$ in $\h_4$, we do not understand its topology. The
situation only gets worse in higher genus.

\section{Preliminaries}

Suppose that $2g-2+n > 0$. Fix a compact oriented surfaced $S$ of genus $g$ 
and a finite subset $P$ of $n$ distinct points in $S$. The corresponding
mapping class group is
$$
\G_{g,n} = \pi_0\Diff^+ (S,P).
$$
By a complex curve, or a curve for short, we shall mean a Riemann surface.
Denote the Teichm\"uller space of marked, $n$-pointed, compact genus $g$ curves
by $\X_{g,n}$. As a set $\X_{g,n}$ is
$$
\left\{
\parbox{2.15in}{orientation preserving diffeomorphisms $f: S \to C$
to a complex curve}
\right\}
\bigg/
\text{isotopies, constant on $P$}.
$$
This is a complex manifold of dimension $3g-3+n$. 

The mapping class group $\G_{g,n}$ acts properly discontinuously on $\X_{g,n}$.
The quotient $\G_{g,n}\bs \X_{g,n}$ is the moduli space $\M_{g,n}$ of
$n$-pointed curves of genus $g$.

Set $H_R = H_1(S;R)$. The intersection pairing $H_R^{\otimes 2} \to R$ is
skew symmetric and unimodular. Set
$$
\Sp(H_R) = \Aut(H_R,\text{intersection pairing}).
$$
The choice of a symplectic basis of $H_R$ gives an isomorphism $\Sp(H_R) \cong
\Sp_g(R)$ of $2g\times 2g$ symplectic matrices with entries in $R$. The action
of $\G_{g,n}$ on $S$ induces a homomorphism
$$
\rho : \G_{g,n} \to \Sp(H_\Z)
$$ 
which is well-known to be surjective.

The {\em Torelli group} $T_{g,n}$ is defined to be the kernel of $\rho$. It is
torsion free. The Torelli space $\T_{g,n}$ is defined by
$$
\T_{g,n} = T_{g,n}\bs \X_{g,n}.
$$ 
It is the moduli space of $n$-pointed Riemann surfaces $(C;x_1,\dots,C_n)$ of
genus $g$ together with a symplectic basis of $H_1(C;\Z)$. Since $T_{g,n}$ is
torsion free, it acts freely on Teichm\"uller space. Torelli space $\T_{g,n}$
is thus a model of the classifying space $BT_{g,n}$. Note that $T_{1,1}$ is
trivial and that $\X_{1,1} = \T_{1,1}$ is just the upper half plane, $\h_1$.
Torelli spaces also exist in genus 0 provided $n\ge 3$. In this case, $\T_{0,n}
= \M_{0,n}$.

The Siegel upper half space of rank $g$
$$
\h_g := \Sp_g(\R)/U(g)
\cong \{\Omega \in \M_g(\C) : \Omega = \Omega^T, \Im \Omega > 0\}
$$
is the symmetric space of $\Sp_g(\R)$. It is a complex manifold of dimension
$g(g+1)/2$. It can usefully be regarded as the moduli space of $g$-dimensional
principally polarized abelian varieties $(A,\theta)$ together with a symplectic
(with respect to the polarization $\theta : H_1(A;\Z)^{\otimes 2}\to \Z$) basis
of $H_1(A;\Z)$. The group $\Sp_g(\Z)$ acts on the framings, and the quotient
$\Sp_g(\Z)\bs \h_g$ is the moduli space $\A_g$ of principally polarized abelian
varieties of dimension $g$.

The decoration $n$ in $\G_{g,n}$, $T_{g,n}$, $\T_{g,n}$, etc will be omitted
when it is zero.

\section{Geography}
\label{geography}

\subsection{The period mapping}
A {\em framed} genus $g$ curve is a compact Riemann surface of genus $g$
together with a symplectic basis $(a_1,\dots,a_g,b_1,\dots,b_g)$ of
$H_1(C;\Z)$. For each framed curve, there is a unique basis $w_1,\dots,w_g$ of
the holomorphic differentials $H^0(C,\Omega^1_C)$ on $C$ such that
$$
\int_{a_j} w_k = \delta_{j,k}.
$$
The period matrix of $(C;a_1,\dots,a_g,b_1,\dots,b_g)$ is the $g\times g$
matrix
$$
\Omega = \left(\int_{b_j} w_k\right).
$$
It is symmetric and has positive definite imaginary part.

As remarked above, $\T_g$ is the moduli space of framed genus $g$
curves. The {\em period mapping}
$$
\T_g \to \h_g
$$
takes a framed curve to its period matrix. It is holomorphic and descends to
the mapping
$$
\M_g = \Sp_g(\Z)\bs \T_g \to \Sp_g(\Z)\bs \h_g = \A_g
$$
that takes the point $[C]$ of $\M_g$ corresponding to a curve $C$ to the 
point $[\Jac C]$ of $\A_g$ corresponding to its jacobian.

\subsection{The jacobian locus} This is defined to be the image $\J_g$ of the
period mapping $\T_g \to \h_g$. It is a locally closed subvariety of $\h_g$. It
is important to note, however, that it is {\em not} closed in $\h_g$.

To explain this, we need to introduce the locus of {\em reducible principally
polarized abelian varieties}. This is
$$
\A_g^\red = \bigcup_{h=1}^{\lfloor g/2 \rfloor}
\im \{\mu_h : \A_h \times \A_{g-h} \to \A_g\}
$$
where $\mu_h(A,B) = A\times B$.\footnote{Note that $\im \mu_h = \im \mu_{g-h}$.
This is why we need only those $h$ between $1$ and $\lfloor g/2 \rfloor$.} Set
$$
\h_g^\red = \text{inverse image in $\h_g$ of } \A_g^\red =
\bigcup_{h=0}^{\lfloor g/2 \rfloor}
\bigcup_{\phi\in \Sp_g(\Z)/\Sp_g(\Z)_{\h_h\times \h_{g-h}}}
\phi \big(\h_h \times \h_{g-h}\big)
$$
where $\h_h \times \h_{g-h}$ is imbedded into $\h_g$ by taking
$(\Omega_1,\Omega_2)$ to the matrix $\Omega_1 \oplus \Omega_2$; the group
$Sp_g(\Z)_{\h_h\times \h_{g-h}}$ denotes the subgroup of $\Sp_g(\Z)$ that fixes
it set-wise.

For a subset $N$ of $\h_g$, set $N^\red = N \cap \h_g^\red$. By
\cite[Prop.~6]{hain:sp3},
$$
\J_g = \J_g^c - \J_g^{c,\red}
$$
from which it follows that
$$
\J_g^{c,\red} =
\bigcup_{\phi \in \Sp_g(\Z)} \bigcup_{h=1}^{\lfloor g/2 \rfloor}
\phi \big(\J_h^{c,\red}\times \J_{g-h}^{c,\red}\big).
$$

\subsection{Curves of compact type}
A genus $g$ curve $C$ of {\em compact type} is a connected, compact, nodal
curve\footnote{A nodal curve is a complex analytic curve, all of whose
singularities are nodes --- that is, of the form $zw = t$.} satisfying:
\begin{enumerate}

\item the dual graph of $C$ is a tree --- this guarantees that the jacobian
$\Jac C$ of $C$ is a principally polarized abelian variety;

\item the sum of the genera of the components of $C$ is $g$.

\end{enumerate}
These are precisely the stable curves of genus $g$ whose generalized jacobian
is compact. The generalized jacobian a genus $g$ curve of compact type is the 
product of the jacobians of its components. It is an abelian variety of
dimension $g$.

An $n$-pointed nodal curve of genus $g$ is a nodal genus $g$ curve $C$ together
with $n$ labeled points in the smooth locus of $C$. An $n$-pointed nodal curve
$(C,P)$ of genus $g$ is {\em stable} if its automorphism group is finite. This
is equivalent to the condition that each connected component of
$$
C - (C^\sing \cup P)
$$
has negative Euler characteristic.

Using the deformation theory of stable curves, one can enlarge $\T_{g,n}$ to
the moduli space $\T_{g,n}^c$ of framed stable $n$-pointed genus $g$ curves of
compact type. This is a complex manifold that contains $\T_g$ as a dense open
subset and on which $\Sp_g(\Z)$ acts (via its action on framings). The quotient
$\Sp_g(\Z)\bs \T_{g,n}^c$ is the space $\M_{g,n}^c$ of stable $n$-pointed,
genus $g$ curves of compact type.

Note that when $n\ge 3$, $\T_{0,n}^c = \Mbar_{0,n}$, the moduli space of stable
$n$-pointed curves of genus 0.

\begin{proposition}
If $2g-2+n>0$, then $\T_{g,n}^c$ is a smooth complex analytic variety of
complex dimension $3g-3+n$ and
$$
\T_{g,n} = \T_{g,n}^c - Z
$$
where $Z$ is a countable union of transversally intersecting smooth divisors.
\end{proposition}

The strata of $Z$ of complex codimension $k$ are indexed by the $k$-simplices
of the quotient $T_{g,n}\bs K^\sep(S,P)$ of the complex of separating curves
$K^\sep(S,P)$ of the $n$-pointed reference surface $(S,P)$ by the Torelli
group. This correspondence will be made more explicit in the next paragraph.

Since $\T_{g,n}=\T_{g,n}^c-Z$, where $Z$ has complex codimension 1, the mapping
$$
\pi_1(\T_{g,n},\ast) \to \pi_1(\T_{g,n}^c,\ast)
$$
is surjective. Its kernel is generated by the conjugacy classes of small loops
about each of the components of $Z$. But these are precisely the conjugacy
classes of Dehn twists on separating simple closed curves (SCCs).

\begin{proposition}
For all $(g,n)$ satisfying $2g-2+n>0$,
$$
\pi_1(\T_{g,n}^c,\ast) \cong
T_{g,n}/\{\text{subgroup generated by Dehn twists on separating SCCs}\}.
$$
\end{proposition}

As in the previous section, we set $H_\Z = H_1(S;\Z)$, where $S$ is the genus
$g$ reference surface. Denote the image of $u\in \Lambda^3 H_\Z$ under the
quotient mapping
$$
\Lambda^3 H_\Z \to \big(\Lambda^3 H_\Z\big)/\big(\theta \wedge H_\Z\big)
$$
by $\ubar$.

The surjectivity of the ``Johnson homomorphism'' $\tau : T_{g,1} \to \Lambda^3
H_\Z$ and Johnson's result \cite{johnson:scc} that its kernel is generated by
Dehn twists on separating SCCs implies quite directly that:

\begin{corollary}
If $g \ge 1$ and $2g-2+n>0$, then
$$
\pi_1(\T_{g,n}^c,\ast) = H_1(\T_{g,n}^c;\Z) \cong
\big\{
(u_1,\dots,u_n) \in \big(\Lambda^3 H_\Z \big)^n : \ubar_1 = \dots = \ubar_n
\big\},
$$
which is a torsion free abelian group.
\end{corollary}

\subsection{The complex of {separating} curves, $K^\sep(S,P)$}

As above $S$ is a compact oriented surface of genus $g$ and $P$ is a subset of
cardinality $n$, where $2g-2+n>0$. A simple closed curve $\gamma$ in $S-P$ is 
separating if $S-(P\cup \gamma)$ is not connected. An SCC is {\em cuspidal} if
it bounds a disk in $S$.

The simplicial complex $K^\sep(S,P)$ has vertices the isotopy classes of
separating SCCs $\gamma$ in $S-P$ that are not cuspidal. The isotopy classes of
the SCCs $\gamma_0,\dots,\gamma_k$ of non-cuspidal separating SCCs span a
$k$-simplex of $K^\sep(S,P)$ if they are disjoint and lie in distinct isotopy
classes. When $P$ is empty, we shall abbreviate $K^\sep(S,P)$ by $K^\sep(S)$.

The correspondence between $T_{g,n}$ orbits of $K^\sep(S,P)$ and the strata of
$\T_{g,n}^{c,\red}$ is given as follows. Given a $k$-simplex $\bgamma =
(\gamma_0,\dots,\gamma_k)$ of $K^\sep(S,P)$, one can contract each of the SCCs
$\gamma_j$. The resulting space $(S/\bgamma,P)$ is the topological model of a
stable, $n$-pointed, genus $g$ (complex) curve of compact type. Every
topological type of stable complex curve of compact type arises in this way.

A marked $n$-pointed, genus $g$ curve of compact type is a homotopy class of
homeomorphisms
$$
(S/\bgamma,P) \to (C,\{x_1,\dots,x_n\})
$$
to an $n$-pointed, genus $g$ stable curve $(C,\{x_1,\dots,x_n\})$. For each
$\bgamma \in K^\sep(S,P)$, one can add a connected ``rational boundary
component'' $\X_\gamma$ to the Teichm\"uller space $\X_{g,n}$ of $(S,P)$ to
obtain a topological space . The mapping class group $\G_{g,n}$ acts on this
enlarged Teichm\"uller space $\X_{g,n}^c$ and the quotient is $\M_g^c$ from
which it follows that $\T_{g,n}^c = T_{g,n}\bs \X_{g,n}^c$.

The stratum of $\T_g^{c,\red}$ that has codimension $k$ in $\T_g^c$ is the
locus of stable curves of compact type with precisely $k$ double points. These
correspond to the $k-1$ simplices of $K^\sep(S,P)$. It follows that the  strata
of $\T_{g,n}^c$ correspond to the simplices of $T_{g,n}\bs K^\sep(S,P)$.

\begin{remark}
Farb and Ivanov \cite{farb-ivanov} have shown that $K^\sep(S)$ is connected
whenever $g\ge 3$.
\end{remark}

\subsection{Singularities and dimension}

Note that $\J_g$ is the quotient of $\T_g$ by the involution
$$
\sigma :
(C;a_1,\dots,a_g,b_1,\dots,b_g) \to (C;-a_1,\dots,-a_g,-b_1,\dots,-b_g)
$$
Since
$$
(C;a_1,\dots,a_g,b_1,\dots,b_g) \cong (C;-a_1,\dots,-a_g,-b_1,\dots,-b_g)
$$
if and only if $C$ is hyperelliptic, this mapping is ramified along the locus
$\H_g$ of hyperelliptic curves. Since this has complex codimension $g-2$ in
$\T_g$, we know that $\J_g$ is singular along the locus of hyperelliptic
jacobians when $g\ge 4$. It is the quotient of the manifold $\T_g$ by $\Z/2$,
so $\J_g$ is always a $\Z[1/2]$-homology manifold.

The period mapping $\T_g^c \to \J_g^c$ has positive dimensional fibers over
$\J_g^{c,\red}$ when $g\ge 3$. As a result, $\J_g^c$ is singular along
$\J_g^{c,\red}$ when $g\ge 4$, and is not a rational homology
manifold.\footnote{An explicit description of the links of the singularities of
the top stratum of $\J_g^{c,\red}$ can be found in \cite[Prop.~6.5]{hain:oort}.
There are similar descriptions in higher codimension.}

Since $\J_g^c$ has dimension $3g-3$ and $\h_g$ dimension $g(g+1)/2$, $\J_g^c$
is a proper subvariety of $\h_g$ when $g\ge 4$. This and the fact that
$\J_g^c$ is not a rational homology manifold when $g\ge 4$ help explain the
difficulty of generalizing Mess' arguments to any $g\ge 4$.

On the positive side, we can say that, since $\J_g^c$ is a closed analytic
subvariety of the Stein manifold $\h_g$, it is a Stein space. Consequently, by
a result of Hamm \cite{hamm} we have:

\begin{proposition}
$\J_g^c$ has the homotopy type of a CW-complex of dimension at most $3g-3$.\qed
\end{proposition}

\section{Homological Tools}

This section can be omitted or skimmed on a first reading.

\subsection{A Gysin sequences} Since $\T_g$ is obtained from the manifold
$\T_g^c$ by removing a countable union of closed subvarieties, we have to be a
little more careful than usual when constructing the Gysin sequence.

Suppose that $X$ is a PL manifold of dimension $m$ and that $Y$ is a closed PL
subset of $X$. Suppose that $A$ is any coefficient system. For each compact
(PL) subset $K$ of $X$, we have the long exact sequence
\begin{multline*}
\cdots \to H^{m-k-1}(Y,Y-(Y\cap K);A) \to H^{m-k}(X,Y\cup(X-K);A)
\cr
\to H^{m-k}(X,X-K;A) \to H^{m-k}(Y,Y-(Y\cap K);A) \to \cdots
\end{multline*}
of the triple $(X,Y\cup(X-K),X-K)$. Taking the direct limit over all such $K$,
we obtain the exact sequence
$$
\cdots \to H^{m-k-1}_c(Y;A) \to H^{m-k}_c(X-Y;A)
\to H^{m-k}_c(X;A) \to H^{m-k}_c(Y;A) \to \cdots
$$
When $X$ is oriented, we can apply Poincar\'e duality to obtain the following
version of the Gysin sequence:

\begin{proposition}
\label{gysin}
If $A$ is any coefficient system, $X$ an oriented PL-manifold of dimension
$m$ and $Y$ is a closed PL subset of $X$, then there is a long exact sequence 
$$
\cdots \to H^{m-k-1}_c(Y;A) \to H_k(X-Y;A)
\to H_k(X;A) \to H^{m-k}_c(Y;A) \to \cdots
$$
\end{proposition}

\subsection{A spectral sequence}
In practice, we are also faced with the problem of computing $H^\dot_c(Y)$ in
the Gysin sequence above when $Y$ is singular. Suppose that
$$
Y = \bigcup_{\alpha \in I} Y_\alpha
$$
is a locally finite union of closed PL subspaces of a PL manifold $X$. Set
$$
Y_{(\alpha_0,\dots,\alpha_k)} :=
Y_{\alpha_0} \cap Y_{\alpha_1} \cap \dots \cap Y_{\alpha_k}
$$
and
$$
\Y_k = \coprod_{\alpha_0 < \alpha_1 < \dots < \alpha_k}
Y_{(\alpha_0,\dots,\alpha_k)}.
$$
The inclusions
$$
d_j :  Y_{(\alpha_0,\dots,\alpha_k)} \hookrightarrow
Y_{(\alpha_0,\dots, \widehat{\alpha_j},\dots,\alpha_k)}
$$
induce face maps
$$
d_j : \Y_k \to \Y_{k-1} \quad j = 0, \dots ,k
$$
With these, $\Y_\dot$ is a strict simplicial space.

\begin{proposition}
There is a spectral sequence
$$
E_1^{s,t} \cong H^t_c(\Y_s;A) \implies H^{s+t}_c(Y;A)
$$
whose $E_1$ differential
$$
d_1 : H^t_c(\Y_s;A) \to H^t_c(\Y_{s+1};A)
$$
is $\sum_j (-1)^j d_j^\ast$.
\end{proposition}

\begin{proof}
This follows rather directly from the standard fact that the natural
``augmentation''
$$
\epsilon : |\Y_\dot| \to Y
$$
from the geometric realization\footnote{This is the quotient of $\coprod_{k\ge
0} \Y_k \times \Delta^k$ by the natural equivalence relation generated by
$(y,\partial_j\xi) \cong (d_j y, \xi)$, where $(y,\xi) \in \Y_k \times
\Delta^{k-1}$, where $\Delta^d$ denotes the standard $d$-simplex, and
$\partial_j : \Delta^{k-1} \to \Delta^k$ is the inclusion of the $j$ th face.}
is a homotopy equivalence. Since $Y$ is a locally finite union of the closed
subspaces $Y_\alpha$, the natural mapping
$$
S_c^\dot(Y;A) \to S^\dot_c(\Y_\dot;A)
$$
is a quasi-isomorphism. The spectral sequence is that of the double complex
$S^\dot_c(\Y_\dot;A)$. The quasi-isomorphism implies that the spectral sequence
abuts to $H^\dot_c(Y;A)$.
\end{proof}

Suppose now that each $Y_\alpha$ is an oriented PL submanifold of dimension $2d$
and codimension 2 in $X$. In addition, suppose that the components of $Y_\alpha$
of $Y$ intersect transversally in $X$, so that each component of $\Y_s$ has
dimension $2d-2s$. By duality, the $E_1$ term of the spectral sequence can we
written as
$$
E_1^{s,t} = H_{2d-2s-t}(Y_s;A)
$$
The differential
$$
d_1 : H_{2d-2s-t}(\Y_s;A) \to H_{2d-2s-t-2}(\Y_{s+1};A)
$$
is the alternating sum of the Gysin mappings $d_j^\ast : \Y_{s+1} \to \Y_s$.

\begin{proposition}
There is a natural isomorphism $H^{2d}_c(Y;A) \cong H_0(\Y_0;A)$ and an exact
sequence
$$
H^{2d-2}_c(Y;A) \to H_2(\Y_0;A) \to H_0(\Y_1;A) \to H^{2d-1}_c(Y;A)
\to H_1(\Y_0;A) \to 0
$$
\end{proposition}

\subsection{Cohomology at infinity}

For a topological space $X$ and an $R$-module (or local coefficient system of
$R$-modules) $A$, set
$$
H^\dot_\infty(X;A) = \liminj{\substack{K\subseteq X \cr\text{ compact}}}
H^\dot(X-K;A).
$$
We shall call it the {\em cohomology at infinity} of $X$. When $X$ is a compact
manifold with boundary, then there is a natural isomorphism
$$
H^\dot_\infty(X-\partial X;A) \cong H^\dot(\partial X; A).
$$

\begin{remark}
I am not sure if this definition appears in the literature, although I would be
surprised if it does not. Similar ideas have long appeared in topology, such as
in the paper Bestvina and Feighn \cite{bestvina}, where they introduce the
notion of a space being ``$r$-connected at infinity.''
\end{remark}

The direct limit of the long exact sequence of the pairs $(X,X-K)$, where $K$
is compact, is a long exact sequence
\begin{equation}
\label{les}
\cdots \to
H^{k-1}_\infty(X;A) \to H^k_c(X;A) \to H^k(X;A) \to H^k_\infty(X;A)
\to \cdots
\end{equation}

When $X$ is an oriented manifold of dimension $m$, Poincar\'e duality gives an
isomorphism
$$
H_{m-k}(X;A) \to H^k_c(X;A).
$$
Thus, if $X$ has the homotopy type of a CW-complex of dimension $d$, then
$$
H^k_c(X;A) = 0 \text{ when } k < m-d \text{ and }
H^k(X;A) = 0 \text{ when } k>d.
$$
Plugging these into the long exact sequence (\ref{les}), we see that
$$
H^k_\infty(X;A) \cong H^{k}(X;A)  \text{ when } k < m-d-1
$$
and
$$
H^{k}_\infty(X;A) \cong H^{k+1}_c(X;A) \text{ when } k>d.
$$
Moreover, since
$$
H^k(X;R) \to \Hom_R\big(H^{m-k}_c(X;R),R\big)
$$
is an isomorphism, we have an isomorphism
$$
H^k_\infty(X;R) \to \Hom_R\big(H^{m-k-1}_\infty(X;R), R\big)
$$
when $k < m-d-1$.

One can consider vanishing of $H_\infty^\dot(X)$ modulo the Serre class of
finitely generated $R$-modules. This leads us to the following bizarre result.

\begin{lemma}
Suppose that $R$ is a PID. If $H_\infty^k(X;R)$ is a countably generated
$R$-module when $a \le k \le b$, then $H^k(X;R)$ is a finitely generated
$R$-module when $a \le k \le b$.
\end{lemma}

\begin{proof}
Since $X$ has the homotopy type of a countable CW-complex, each $H_k(X;R)$ is
countably generated. Poincar\'e duality then implies that $H_c^k(X;R)$ is also
countably generated. The Universal Coefficient Theorem implies that $H^k(X;R)$
surjects onto $\Hom_R(H_k(X;R),R)$. Consequently, if $H^k(X;R)$ is countably
generated, it is finitely generated. The result now follows from the exact
sequence (\ref{les}).
\end{proof}

\section{Discussion and Problems}

One natural approach to the problem of understanding the topology of $\T_g$ is
to view it as the complement of the normal crossings divisor $\T_g^{c,\red}$ in
$\T_g^c$. The space $\T_g^c$ can in turn be studied via the period mapping
$\T_g^c \to \J_g^c$. There are two ways to factorize this. The first is to take
the quotient $\cQ_g^c:= \T_g^c/\langle \sigma \rangle$ of $\T_g^c$ by the
involution
$$
\sigma : [C;a_1,\dots,a_g,b_1,\dots,b_g]
\to [C;-a_1,\dots,-a_g,-b_1,\dots,-b_g].
$$
The mapping $\T_g^c \to \cQ_g^c$ is branched along the hyperelliptic locus
$\H_g^c$. The second is to consider the Stein factorization (cf.\
\cite[p.~213]{grauert-remmert})
$$
\T_g^c \to \S_g^c \to \J_g^c
$$
of the period mapping. The important properties of this are that $\S_g^c$ is a
complex analytic variety, the first mapping has connected fibers, while the
second is finite (in the sense of analytic geometry). The two factorizations
are related by the diagram
$$
\xymatrix{
\T_g^c \ar[rr] \ar[d] \ar[drr]|{\text{period map}} && \S_g^c \ar[d]\cr
\cQ_g^c \ar[rr] && \J_g^c
}
$$
where all spaces are complex analytic varieties, all mappings are proper,
holomorphic and surjective. The horizontal mappings have connected fibers and
the vertical mappings are finite and two-to-one except along the hyperelliptic
locus.

Perhaps the first natural problem is to understand the topology of $\J_g^c$.

\begin{problem}[{\em Topological Schottky Problem}]
Understand the homotopy type of $\J_g^c$ and use it to compute $H^\dot(\J_g^c)$
and $H_c^\dot(\J_g^c)$.
\end{problem}

The first interesting case is when $g=4$, where $\J_g^c$ is a singular divisor
in $\h_4$. It would also be interesting and natural to compute the intersection
homology of $\J_g^c$.

A knowledge of the topology of $\J_g$, $\J_g^c$ or $\cQ_g^c$ should help with
the computation of the $\sigma$-invariant homology
$H_\dot(\T_g)^{\langle\sigma\rangle}$.\footnote{When $g=3$, $\cQ_3$ is obtained
from $\h_3$ by first blowing up the singular locus of $\h_3^\red$, which is
smooth of complex codimension 3, and then blowing up the proper transforms of
the components of $\h_3^\red$. The strata of $\h_3^\red$ are described in
detail in \cite{hain:sp3}.}  When $2$ is invertible in the coefficient ring
$R$, we can write
$$
H_\dot(\T_g;R) = H_\dot(\T_g;R)^+ \oplus H_\dot(\T_g;R)^-
$$
where $\sigma$ acts as the identity on $H_\dot(\T_g;R)^+$ and as $-1$ on
$H_\dot(\T_g;R)^-$.

\begin{problem}
Determine whether or not $H_\dot(\T_g;\Z[1/2])^-$ is always a finitely
generated $\Z[1/2]$-module. Does the infinite topology of $\T_g$ comes from
$\J_g$?
\end{problem}

To get one's hands on $H_\dot(\T_g)^-$, it is necessary to better understand
the topology of $\T_g^c$ or $\S_g^c$. Since $\J_g^c$ is a Stein space, so is
$\S_g^c$. Hamm's result \cite{hamm} (see also \cite[p.~152]{smt}) implies that
$\S_g^c$ has the homotopy type of a CW-complex of dimension at most $3g-3$.
Since $\S_g^c$ is not a rational homology manifold when $g\ge 3$, it is
probably most useful to study the topology of the manifold $\T_g^c$ as
Poincar\'e duality will then be available.\footnote{If one uses intersection
homology instead, then duality will still be available. For this reason it is
natural to try to compute the intersection homology of $\S_g^c$ and $\J_g^c$.}

\begin{problem}
Determine good bounds for the homological dimension (or the CW-dimension) of
$\T_g^c$.
\end{problem}

The best upper bound that I know of is obtained using Stratified Morse Theory:

\begin{proposition}
If $g\ge 1$ and $2g-2+n>0$, then the dimension of $\T_{g,n}^c$ as a CW-complex
satisfies
$$
2(g-2+n) \le \dim_\CW \T_{g,n}^c \le
\bigg\lfloor \frac{7g-8}{2}\bigg\rfloor + 2n.
$$
\end{proposition}

\begin{proof}
Suppose that $g\ge 2$. Since $\T_{g,n}^c \to \T_g^c$ is proper with fibers of
complex dimension $n$,
$$
\dim_\CW \T_{g,n}^c = 2n + \dim_\CW \T_g^c.
$$
To establish the lower bound, we need to exhibit a topological $(2g-4)$-cycle in
$\T_g^c$ that is non-trivial in $H_\dot(\T_g^c)$. Observe that
$$
Y := (\T_{1,1}^c)^2 \times (\T_{1,2}^c)^{g-2}
$$
is a component of the closure of the locus of chains of $g$ elliptic curves in
$\T_g^c$. Let $E$ be any elliptic curve. Then $E\subset \T_{1,2}$ and $Y$
therefore contains the projective subvariety $E^{g-2}$. The rational homology
class of a compact subvariety of a K\"ahler manifold is always non-trivial
(just integrate the appropriate power of the K\"ahler form over it). Since
$\T_g^c$ covers $\M_g^c$, which is a K\"ahler orbifold, $\T_g^c$ is a K\"ahler
manifold and the class of $E^{g-2}$ is non-trivial in $H_\dot(\T_g^c)$. This
establishes the lower bound.

The upper bound is a direct application of Stratified Morse Theory \cite{smt}.
Since $\J_g^c$ is a Stein space, it is a closed analytic subvariety of $\C^N$
for some $N$. The constant $c$ in \cite[Thm.~1.1*, p.~152]{smt} is thus 1.
Consequently, $\dim_\CW \T_{g,n}$ is bounded by
$$
d(g,n) := \sup_{k\ge 0} \big(2k + f(k)\big)
$$
where $f(k)$ is the maximal complex dimension of a subvariety of $\T_{g,n}^c$
over which the fiber of $\T_{g,n}^c \to \J_g^c$ is $k$ dimensional. Observe that
$d(g,n) = d(g) + 2n$, when $g\ge 2$, where $d(g) = d(g,0)$. Since $\T_{1,n}^c\to
\h_1$ is proper with fibers of dimension $n-1$, we have $d(1,n)=1 + 2(n-2) =
2n-1$.

Since the mapping $\T_g \to \J_g$ is finite and since the components of
$\T_g^{c,\red} \to \J_g^{c,\red}$ are the $\Sp_g(\Z)$ orbits of the period
mappings
$$
\T_{h,1}^c\times \T_{g-h,1}^c \to \J_h^c\times \J_{g-h}^c
$$
we have
$$
d(g) =
\max\big[3g-3,\max\big\{d(h,1)+d(g-h,1):1\le h \le g/2\big\}\big].
$$
The formula for the upper bound $d(g)$ is now easily proved by induction on $g$.
\end{proof}

The upper bound is not sharp when $g\le 2$; for example, $\T_{1,1}^c=\h_1$ and
$\T_2^c = \h_2$, which are contractible. I suspect it is not sharp in higher
genus as well. The first interesting case is to determine the CW-dimension of
$\T_3^c$.

The upper bound on the homological dimension of $\T_g^c$ implies that
$$
H_k(\T_g^c) \cong H_\infty^{6g-7-k}(\T_g^c) \text{ when }
k < \bigg\lceil \frac{5g-6}{2}\bigg\rceil
$$
So the low dimensional homology of $\T_g^c$ is related to the topology at
infinity of $\T_g^c$.

\begin{problem}
Try to understand the ``topology at infinity'' of $\T_g^c$. In particular, try
to compute $H_\infty^k(\T_g^c)$ for $k$ in some range $k\ge d_o$.
Alternatively, try to compute the homology $H_\dot(\T_g^c)$ in lower
degrees.
\end{problem}

Note that $\T_2^c$ is a manifold with boundary $S^5$. The boundary of $\J_3^c$
is $S^{11}$. The first interesting case is in genus 3.

\begin{problem}
Compute $H^\dot_\infty(\T_3^c)$.
\end{problem}

The homology of $\T_g^c$ is related to that of $\T_g$ via the Gysin sequence.
In order to apply it, one needs to understand the topology of the divisor
$\T_g^{c,\red}$. This is built up out of products of lower genus Torelli spaces
of compact type. The combinatorics of the divisor is given by the complex
$K^\sep(S)/T_g$.

\begin{problem}
Compute the $\Sp_g(\Z)$-module $H^k_c(\T_g^{c,\red})$ in some range $k\ge k_o$.
\end{problem}

We already know that $H^{6g-7}_c(\T_g^{c,\red}) = H_0(\B_g^c)$ and that there
is a surjection $H^{6g-8}_c(\T_g^{c,\red})\to H_1(\B_g^c)$, where $\B_g^c$
denotes the normalization (i.e., disjoint union of the irreducible components)
of $\T_g^{c,\red}$. In concrete terms:
$$
\B_g^c = \coprod_{h=1}^{\lfloor g/2 \rfloor}
\coprod_{\phi \in \Sp_g(\Z)_{(\T_{h,1}^c\times \T_{g-h,1}^c)}}
\phi\big(\T_{h,1}^c\times \T_{g-h,1}^c \big).
$$

Lurking in the background is the folk conjecture
$$
H_k(\T_g) \text{\em is finitely generated when } k < g-1.
$$
If true, this places strong conditions on the finiteness of the topology of
$\T_g^c$ and $\T_g^{c,\red}$. It is worthwhile to contemplate (for $g\ge 3$)
the Gysin sequence:
$$
\xymatrix{
\cdots \ar[r]
& H^{6g-10}_c(\T_g^{c,\red}) \ar[r] & H_3(\T_g) \ar[r] & H_3(\T_g^c) \ar[ddll]
\ar@{=}[d]
\cr
&&& H^{6g-10}_c(\T_g^c)
\cr
& H^{6g-9}_c(\T_g^{c,\red})\ar@{->>}[d] \ar[r] & H_2(\T_g) \ar[r] & H_2(\T_g^c)
\ar[ddll] \ar@{=}[d]
\cr
& H_1(\B_g^c) && H^{6g-9}_c(\T_g^c)
\cr
& H^{6g-8}_c(\T_g^{c,\red}) \ar@{=}[d] \ar[r] & H_1(\T_g) \ar[r]^\tau &
H_1(\T_g^c) \ar[r]\ar@{=}[d] & 0
\cr
& H_0(\B_g^c) && \big(\Lambda^3 H\big)/H
}
$$
Here $\tau$ denotes the Johnson homomorphism, realized as the map on $H_1$
induced by the inclusion $\T_g \hookrightarrow \T_g^c$.

Finally, it is interesting to study the topology of the branching locus of
$\T_g^c \to \J_g^c$. This is the locus $\H_g^c$ of hyperelliptic curves of
compact type. Using A'Campo's result \cite{acampo} that the image of the
hyperelliptic mapping class group\footnote{The hyperelliptic mapping class is
the centralizer of a hyperelliptic involution in $\G_g$. It is the orbifold
fundamental group of the moduli space of smooth hyperelliptic curves of genus
$g$.} $\Delta_g$ in $\Sp_g(\Z)$ contains the level two subgroup
$$
\Sp_g(\Z)[2] = \{ A \in \Sp_g(\Z)[2] : A \equiv I \bmod 2\}
$$
of $\Sp_g(\Z)$ and the fact that the image $\Delta_g$ in $\Sp_g(\F_2)$ is
$S_{2g-2}$, the symmetric group on the Weierstrass points, one can see that
$\H_g^c$ has
$$
|\Sp_g(\F_2)|/|S_{2g+2}| = \frac{2^{g^2}\prod_{k=1}^g (2^{2k}-1)}{(2g+2)!}
$$
components. Each component of $\H_g^c$ is smooth and immerses in $\h_g$ via the
period mapping. The irreducible components of $\H_g$ are disjoint in $\J_g$. In
genus 3, $\H_3^c$ has 36 components. Their images are cut out by the 36 even
theta nulls $\vartheta_\alpha : \h_3 \to \C$.

\begin{conjecture}
Each component of $\H_g^c$ is simply connected.
\end{conjecture}

This is trivially true in genus 2, where there is one component which is all of
$\h_2$. If true in genus 3, it implies quite directly the known fact that $T_3$
is generated by $35=36-1$ bounding pair elements. The number $35$ is the rank
of $H_1(\J_3,\H_3)$. The inverse images of generators of $\pi_1(\J_3,\H_3)$,
once oriented, generate $T_3$.\footnote{As the genus increases, the number of
components of $\H_g$ increases exponentially while the minimum number of
generators of $T_g$ increases polynomially. The failure of the genus 3 argument
presented above in higher genus suggests that $\J_g^c$ is not, in general,
simply connected.}

The conjecture has an equivalent statement in more group theoretic terms.
Define the hyperelliptic Torelli group to be the intersection of the
hyperelliptic mapping class group and the Torelli group: 
$$
T\Delta_g := \Delta_g \cap T_g = \ker\{\Delta_g \to \Sp_g(\Z)\}.
$$
It is a subgroup of the Johnson subgroup $K_g := \ker\{\tau : T_g \to \Lambda^3
H/H\}$. Examples of elements in $T\Delta_g$ are Dehn twists on separating
simple closed curves that are invariant under the hyperelliptic involution
$\sigma$. If $\H_{g,\alpha}$ is a component of $\H_g$, then there is an
isomorphism
$$
\pi_1(\H_{g,\alpha},*) \cong T\Delta_g.
$$
The conjugacy classes of twists on a $\sigma$-invariant separating SCC
correspond to loops about components of the divisor $\H_g^{c,\red}$.

Van Kampen's Theorem implies that
$$
\pi_1(\H_{g,\alpha}^c,*) \cong
\pi_1(\H_{g,\alpha})/\text{these conjugacy classes}
\cong T\Delta_g/\text{these conjugacy classes}.
$$
The conjecture is thus equivalent to the statement that $T\Delta_g$ is
generated by the conjugacy classes of Dehn twists on $\sigma$-invariant
separating SCCs.

It is important to understand the topology of the loci of hyperelliptic curves
as it is the branch locus of the period mapping and also because it is
important in its own right.

\begin{problem}
Investigate the topology of $\H_g$ and $\H_g^c$ and their components.
Specifically, compute their homology and the cohomology at infinity of
$\H_{g,\alpha}^c$.
\end{problem}

The period mapping immerses each $\H_{g,\alpha}^c$ in $\h_g$ as a closed
subvariety. Consequently, each $\H_{g,\alpha}^c$ is a Stein manifold and thus
has the homotopy type of a CW-complex of dimension equal to its complex
dimension, which is $2g-1$.

\end{document}